


\magnification=\magstep 1
\parindent = 0 pt
\baselineskip = 16 pt
\parskip = \the\baselineskip

\font\AMSBoldBlackboard = msbm10

\def\RR{{\hbox{\AMSBoldBlackboard R}}}
\def\CC{{\hbox{\AMSBoldBlackboard C}}}

\def\QQ{{\hbox{\AMSBoldBlackboard Q}}}

\def\ZZ{{\hbox{\AMSBoldBlackboard Z}}}

\settabs 12\columns
\rightline{math.NT/9811040}
\vskip 1 true in
{\bf \centerline{SPECTRAL ANALYSIS OF THE LOCAL CONDUCTOR OPERATOR}}
\vskip 0.5 true in
\centerline{Jean-Fran\c{c}ois Burnol}
\par
\centerline{November 1998}
\par
The conductor operator acts on a function through multiplying it with the logarithm of the norm of the variable both in position and in momentum space and adding the outcomes. It makes sense at each completion of an arbitrary number field and arose in previous papers by the author where it was shown to be intimately connected with the Explicit Formula of Analytic Number Theory. I complete here its spectral analysis: the conclusion is that the local contribution to the Explicit Formula expressed as an integral over the critical line is completely equivalent to, indeed precisely realizes this spectral analysis. In this picture the inversion is related to complex conjugation, the additive Fourier transform is closely related to the Tate-Gel'fand-Graev Gamma function on the critical line, and the conductor operator itself is just the logarithmic derivative of the Gamma function, again on the critical line, acting as a multiplier. From this one gets that the conductor operator is not only dilation and Fourier invariant, which was obvious, but also inversion invariant. In this manner, an operator theoretic and spectral interpretation of each local contribution to the Explicit Formula on the critical line has been obtained.\par
\vfill
{\parskip = 0 pt
62 rue Albert Joly\par
F-78000 Versailles\par
France\par}

jf.burnol@dial.oleane.com

\eject

{\bf
TABLE OF CONTENTS\par
\par
The conductor operator and the Explicit Formula\par
The discrete spectrum of the $p-$adic conductor operator\par
The continuous spectrum of the $p-$adic conductor operator\par
The Poisson summation formula from a local functional equation\par
Intrinsic spectral analysis in the general case\par
References\par}
\par
\vfill \eject

{\bf The conductor operator and the Explicit Formula\par}

Ever since Tate's Thesis [1] there has always been a tension between the additive and multiplicative structures in the adelic approach to the Riemann Zeta function and its (abelian) generalizations. Tate gave a proof of the Functional Equation which reduced it to a statement about the additive Fourier transform of a multiplicatively homogeneous distribution. His proof uses the (adelic) Poisson summation formula as a tool. Although one can in turn deduce the classical Poisson formula from the Functional Equation it would be misleading to identify it as the true meaning of the Functional Equation.

In a paper dating back to about the same period Weil [2] showed that the contributions to one side of the Explicit Formulae of Analytic Number Theory could be indexed in terms of the various completions $\nu$ of the number field $K$ and furthermore that they truly had their natural home on the multiplicative groups $K_\nu^\times$. But his ``finite part'' symbol left the comparison between archimedean and finite primes in an unsatisfactory state. And the ramified primes also required a tailor-made derivation.

It is fair to say that close to fifty years after the appearance of these works, their proper unification is still a job for the future, with the Riemann Hypothesis at stake.

A great step forward was accomplished by Haran [3] who, for the Riemann Zeta function, was able to give to all Weil local terms an {\it exactly identical} formulation. He discovered that although the analytical origin of these local terms give them at their birth the shape of a {\it multiplicative} convolution, they also exist as an {\it additive} convolution. In [4] I have shown how to prove Haran's Theorem for the general Dirichlet-Hecke $L$-series in a direct manner which does not necessitate to differentiate among finite and infinite, ramified and non-ramified primes. In the course of this derivation one encounters an interesting operator which I have called the conductor operator, because it computes the conductor of a ramified character when it is applied to a test-function with this homogeneity. In this paper I give the complete spectral analysis of this operator. In a way the job has already been done in [4] and the only problem here is to give proper names to the various elements of the Explicit Formula. Indeed I will show that the local term of the Explicit Formula at the place $\nu$ actually {\it is} the spectral analysis of the conductor operator on $K_\nu$. This only works when the Explicit Formula is expressed as an integral over the critical line $Re(s) = {1\over2}$, which thus appears in a natural manner here.

So let $K$ be a number field and $\nu$ one of its places, either finite or archimedean, with completion $K_\nu$. The conductor operator $H$ acts as an unbounded operator in $L^2(K_\nu, dt)$, but I will leave its exact domain of definition unspecified for the moment. For the time being we will consider its action on the space ${\cal S}_0(K_\nu)$ of Schwartz-Bruhat $\varphi(t)$ that are identically $0$ in a neighborhood of the origin and of compact support. We take $dt$ to be the additive Haar measure which is self-dual for the duality corresponding to the additive character $\lambda(t) = \exp(2\pi i\{Tr_{K_\nu/\QQ_p}(t)\})$ at a finite place, $\lambda(t) = \exp(-2\pi i t)$ at a real place, $\lambda(z) = \exp(-2\pi i (z+\overline{z}))$ at a complex place. The module $|x|$ of $x \in K_\nu$ is defined so that $d(xt) = |x|dt, x\neq 0, |0| = 0$. It can be computed as $|Norm_{K_\nu/\QQ_p}(x)|_p$ at a finite place, with $|p|_p = {1\over p}$, as the usual absolute value on $\RR$, as the square of the usual absolute value on $\CC$. At a finite place the module of any uniformizer of the ring of local integers is $1/q$ with $q$ being the cardinality of the residue field (which is a power of the prime number $p$ over which $\nu$ lies). One also denotes by $\delta$ the differental exponent at the finite place $\nu$ which is such that the volume of the ring of integers is $q^{-\delta/2}$. The Fourier transform is taken with respect to the basic additive character chosen above.

The conductor operator is defined as follows:
$$H(\varphi)(t) = \log(|t|)\,\varphi(t) + {\cal F}\left(\log(|\xi|)\,{\cal F}^{-1}(\varphi)(\xi)\right)(t)$$
The right-hand side is at any rate a continuous function which belongs to $L^2(K_\nu, dt)$. It is a simple exercise to check:

{\bf Theorem:} $H$ commutes with the Fourier transform and with the isometric action of the multiplicative group $\varphi(t) \mapsto R_x(\varphi)(t) = |x|^{-1/2} \varphi(t/x)$.

There is also the Inversion $I: \varphi \mapsto (t \mapsto {1\over |t|}\varphi({1\over t}))$ which acts on the domain and (isometrically) on the target.

{\bf Theorem:} $H$ commutes with the inversion $I$.

This is less obvious, but will appear so once the spectral analysis of $H$ is completed. The main tool for this analysis is the explicit formula of [4], which is a straightforward derivation from the fundamental relation between homogeneous distributions ([1], [5], [6]):
$${\cal F}(\chi(x)|x|^{s-1}) = \Gamma(\chi, s)\chi^{-1}(y)|y|^{-s}$$
Here $\chi$ is a multiplicative (quasi-)character and $\chi(x)|x|^{s-1}$ and $\chi^{-1}(y)|y|^{-s}$ are considered as (tempered) distributions on $K_\nu$. The function $\Gamma(\chi, s)$ of the character $\chi$ is meromorphic with respect to the complex parameter $s$ and of course depends on them only through their combination $\chi(x)|x|^s$. It is analytic and non-vanishing in the critical strip, and is periodic with period $2\pi i\over \log(q)$ if $\nu$ is finite. It also has the following properties
$$\Gamma(\chi, s)\cdot\Gamma(\chi^{-1}, 1-s) = \chi(-1)$$
$$\overline{\Gamma(\chi, s)} = \chi(-1)\Gamma(\overline{\chi}, \overline{s})$$
The related function $\Lambda(\chi,s) = - {\Gamma^\prime(\chi,s)\over\Gamma(\chi,s)}$ satisfies
$$\Lambda(\chi,s) = \Lambda(\chi^{-1}, 1-s)$$
$$\overline{\Lambda(\chi,s)} = \Lambda(\overline{\chi},\overline{s})$$
Assuming from now on that $\chi$ is a unitary character we get
$$Re(s) = {1\over2}\Rightarrow\left|\Gamma(\chi,s)\right|^2 = 1$$
$$Re(s) = {1\over2}\Rightarrow \Lambda(\chi,s) = \Lambda(\overline{\chi},\overline{s})\hbox{ and is real}$$
For $\chi$ a ramified character with conductor exponent $c$ Tate [1] computed its Gamma function to be simply $q^{(c+\delta)s}$ up to a multiplicative constant whose exact value doesn't concern us here. So, for a ramified character
$$\Lambda(\chi,s) = -(c+\delta)\log(q)$$
For a smooth function $g$ with compact support on $(0,\infty)$ its Mellin Transform $\widehat{g}(s)$ is defined as
$$\widehat{g}(s) = \int_{0}^{\infty} g(u){u}^s \, { du \over u}$$
Starting with such a $g$ and a unitary character $\chi$ and defining $\varphi(t) = g(|t|) \chi^{-1}(t)$ for $t\neq0$, $\varphi(0) = 0$, the explicit formula [4] now stands as $(t\neq0)$
$$H(\varphi)(t) = -{1 \over {2\pi i}}\int_{Re(s)=1/2} \widehat{g}(s)\Lambda(\chi, s)|t|^{-s}\chi^{-1}(t)\, ds$$

{\bf The discrete spectrum of the $p$-adic conductor operator}

For concreteness I will first elucidate the spectrum of the conductor operator over $\QQ_p$. As it commutes with averaging over the units one can split the discussion between the ramified and the non-ramified spectrum. It will turn out that this is also the decomposition between the discrete and the continuous parts of the spectrum.

First the ramified spectrum: from what precedes
$$H(\varphi)(t) = c\log(p)\,\varphi(t)$$
if $\varphi$ has the homogeneity of a ramified character with conductor exponent $c$. Hence:

{\bf Theorem:} The ramified part of the spectrum is purely discrete. It has infinite multiplicity and support $\{\log(p), 2\log(p), 3\log(p), \dots\}$ if $p>2$, $\{2\log(2), 3\log(2), \dots\}$ if $p=2$.

{\bf The continuous spectrum of the $p$-adic conductor operator}

We can now restrict $H$ to the invariant part $L^2_{inv}(\QQ_p)$ spanned by functions depending only on the norm $|t|$. Let's first evaluate $H(\theta)$ for $\theta$ the characteristic function of the units. For this one can use the Fourier transform $G$ of $-\log(|x|)$ which was determined by Vladimirov [7] (see also [4])
$$G(\varphi) = {\log(p) \over 1 - 1/p}\left(\int_{|t|\leq 1}(\varphi(t) - \varphi(0))\, {dt\over|t|} + \int_{|t|>1}\varphi(t)\, {dt\over|t|} + {1\over p}\varphi(0)\right)$$
and the outcome for $H(\theta)(t) = \log(|t|)\theta(t) - (G*\theta)(t)$ is
$$\eqalign{
(|t| < 1)\quad\ H(\theta)(t) &= -\log(p) \cr
(|t| = 1)\quad\ H(\theta)(t) &= \;0 \cr
(|t| > 1)\quad\ H(\theta)(t) &= -{\log(p)\over|t|} }$$
By dilation invariance one has $H(\theta_i)(t) = H(\theta)(p^it)$ for $\theta_i(t) = \theta(p^it)$. So the matrix elements of $H$ in the orthogonal basis $\{\theta_i, i\in\ZZ\}$ depend only on $i-j$. An orthonormal basis of $L^2_{inv}(\QQ_p)$ is given by the functions
$$\eta_i(t) = (1 - 1/p)^{-1/2} p^{-i/2} \theta(p^i t)$$
and the matrix elements of $H$ in this basis are
$$\eqalign{
(i\neq j)\quad\ H_{ij} &=\ <H\eta_i|\eta_j>\ = -\log(p)\,p^{-|i-j|/2} \cr
(i = j)\quad\ H_{ii} &=\ <H\eta_i|\eta_i>\ = 0 }$$
We can now use the isometry from $L^2_{inv}(\QQ_p)$ to $L^2(S^1, {d\theta\over2\pi})$ which sends $\eta_i$ to $e_i(z) : z\mapsto z^i$. Under this isometry one deduces from the above results that $H$ is transported to the multiplication operator with the function $-2\log(p)\ Re({z\over\sqrt{p}-z}) = -\log(p) \sum_{i\neq0}p^{-|i|/2}e_i(z)$. The image of the unit circle under the homographic transformation $z\mapsto {z\over\sqrt{p}-z}$ is again a circle (with center at $1/(p-1)$, radius $\sqrt{p}/(p-1)$) so that:

{\bf Theorem:} The conductor operator is bounded on $L^2_{inv}(\QQ_p)$ and has a purely continuous spectrum with support $[{-2\log(p)\over\sqrt{p}-1}, {+2\log(p)\over\sqrt{p}+1}]$.

The inversion simply exchanges $\eta_i$ with $\eta_{-i}$ so that it is transported to the isometry of $L^2(S^1)$ given by $f(z)\mapsto f(\overline{z})$. Hence:

{\bf Theorem:} $H$ commutes with the inversion $I$.

The proof of the following statement is left to the reader.

{\bf Theorem:} The Fourier transform on $L^2_{inv}(\QQ_p)$ is mapped to the unitary operator on $L^2(S^1)$ given by
$$f(z)\mapsto {\sqrt{p} - z\over\sqrt{p} - \overline{z}}f(\overline{z})$$

{\bf Hint:} one starts from ${\cal F}(\theta) = -{1\over p}\theta_1 + (1 - {1\over p})\sum_{i\leq0}\theta_i$.

Before reworking the previous paragraphs into more intrinsic terms, I will make an apart\'e on one approach, certainly very well-known, to the Poisson summation formula.

{\bf The Poisson summation formula from a local functional equation\par}

Let $x >1$, $L_x(s) = 1/(1-x^{-s})$, $\Lambda_x(s) = - L^\prime_x(s)/L_x(s) = \log(x)\,x^{-s}/(1-x^{-s})$.
$$Re(s) > 0\;\Rightarrow\;\Lambda_x(s) = \log(x)\,\sum_{k\geq1}x^{-ks}$$
$$\Lambda_x(s) + \Lambda_x(-s) = -\log(x)$$
For a smooth function $g$ with compact support on $(0,\infty)$ the first of these equations and Mellin inversion implies
$${1\over{2\pi i}}\int_{Re(s)=1} \widehat{g}(s)\Lambda_x(s)\, ds = \log(x)\,\sum_{k>0}g(x^k)$$
while the second equation (functional equation) gives
$$\eqalign{
{1\over{2\pi i}}\int_{Re(s)=-1} \widehat{g}(s)\Lambda_x(s)\, ds &= -{1\over{2\pi i}}\int_{Re(s)=1} \widehat{g}(-s)(\Lambda_x(s)+\log(x))\, ds \cr
&= -\log(x)g(1) - \log(x)\,\sum_{k<0}g(x^k)
}$$
The poles of $\Lambda_x(s)$ are at $j\,{2\pi i\over\log(x)}, j\in\ZZ$ with residue $+1$ so that one obtains in the end
$$\log(x)\,\sum_{k\in\ZZ}g(x^k) = \sum_{j\in\ZZ}\widehat{g}(j\,{2\pi i\over\log(x)})$$
which is one formulation of Poisson summation formula.

{\bf Intrinsic spectral analysis in the general case\par}

We go back to the general number field $K$ and its completion $K_\nu$ (archimedean or finite). We choose the multiplicative Haar measure $d^*x$ on $K_\nu^\times$ to assign a volume of $1$ to the units if $\nu$ is finite, and to be mapped to $du/u$ under $x\mapsto u = |x|$ for an archimedean place. For the positive constant $\alpha$ chosen such that $d^*x = \alpha^{-2}\,{dx\over |x|}$ the maps
$$\varphi(t)\mapsto\{x\mapsto\alpha\sqrt{|x|}\varphi(x)\}$$
$$f(x)\mapsto\{t\mapsto\alpha^{-1}|t|^{-1/2}f(t)\}$$
establish a unitary isomorphism between $L^2(K_\nu,dt)$ and $L^2(K_\nu^\times,d^*x)$.

If we transport the conductor operator from $L^2(K_\nu,dt)$ to $L^2(K_\nu^\times,d^*x)$ we find that its dilation invariance becomes a commutation property with the multiplicative translations on $K_\nu^\times$. Hence, by the well-known generalization of Fourier Theory, $H$ can be completely analyzed in terms of the (unitary) characters of that locally compact abelian group. Although it is naturally defined in additive terms, its analysis proceeds in multiplicative terms.

So let $S_\nu$ be the dual of $K_\nu^\times$, that is the locally compact abelian group of unitary multiplicative characters. For a finite place $S_\nu$ is the union of countably many circles indexed by the characters of the compact subgroup of units, for the real place it is the union of two real lines indexed by the trivial and the sign characters, for the complex place it is the union of countably many real lines indexed by the characters of the unit circle.

There is a unique Haar measure $\mu$ on $S_\nu$ so that the inversion formula of Fourier Theory and then the Plancherel identity hold:
$$\eqalign{
f(\chi) &= \int_{K_\nu^\times}f(x)\chi(x)\,d^*x \cr
f(x) &= \int_{S_\nu}f(\chi)\,\overline{\chi(x)}\,d\mu(\chi) \cr
\int_{K_\nu^\times}|f(x)|^2\,d^*x &= \int_{S_\nu}|f(\chi)|^2\,d\mu(\chi)
}$$

Under these isometries the Inversion $\varphi(t)\mapsto{1\over|t|}\varphi({1\over t})$ on $L^2(K_\nu,dt)$ ($f(x)\mapsto f(1/x)$ on $L^2(K_\nu^\times,d^*x)$) becomes the operation of complex conjugation $$f(\chi) \mapsto f(\overline{\chi})$$
It is also interesting to look at the fate of the additive Fourier Transform ${\cal F}$. First transported to $L^2(K_\nu^\times,d^*x)$ it maps the function $f(x)$ (chosen ``smooth'' with compact support) to the function $F(x) = \sqrt{|x|}\,{\cal F}(|x|^{-1/2}f(x))$. This implies:
$$F(\chi) = \int_{K_\nu^\times} F(x)\chi(x)\,d^*x = \alpha^{-2}\ \int_{K_\nu} {\cal F}(|t|^{-1/2}f(t))\,\chi(t)|t|^{-1/2}\,dt $$
$$F(\chi) = \alpha^{-2}\ \Gamma(\chi,1/2)\ \int_{K_\nu} |t|^{-1/2}f(t)\,\chi(1/t)|t|^{-1/2}\,dt$$
$$F(\chi) = \Gamma(\chi,1/2)\ f(\overline{\chi})$$

We have obtained the first statements of:

{\bf Theorem:} Under the isomorphims $L^2(K_\nu,dt) \sim L^2(K_\nu^\times,d^*x) \sim L^2(S_\nu,d\mu)$ the inversion is transported to the unitary $$f(\chi)\mapsto f(\overline{\chi})$$ and the additive Fourier transform is transported to the unitary
$$f(\chi) \mapsto F(\chi) = \Gamma(\chi,1/2)\ f(\overline{\chi})$$ while the conductor operator is transported to the hermitian operator
$$f(\chi) \mapsto F(\chi) = - \Lambda(\chi,1/2)\ f(\chi)$$

{\bf Corollary:} The conductor operator commutes with the inversion.

{\bf Proof (of Corollary):} $\Lambda(\chi,1/2) = \Lambda(\overline{\chi},1/2)$.

{\bf Proof (of Theorem):} We start with the explicit formula of [4]
$$H(\varphi)(t) = -{1 \over {2\pi i}}\int_{Re(s)=1/2} \widehat{g}(s)\Lambda(\chi, s)|t|^{-s}\chi^{-1}(t)\, ds$$
where $\chi$ is a unitary character, $g$ is smooth with compact support on $(0,\infty)$, and $\varphi(t) = g(|t|)\chi^{-1}(t)$. Let's write now $f(x) = \sqrt{|x|} g(|x|)\chi^{-1}(x)$ so that the conductor operator transported to $K_\nu^\times$ sends $f$ to
$$F(x) = -{1 \over {2\pi i}}\int_{Re(s)=1/2} \widehat{g}(s)\Lambda(\chi, s)|x|^{{1\over2}-s}\chi^{-1}(x)\, ds$$
$$F(x) = -{1 \over {2\pi}}\int_{-\infty}^{+\infty} \widehat{g}({1\over2} + i\tau)\Lambda(\chi, {1\over2} + i\tau)|x|^{-i\tau}\chi^{-1}(x)\, d\tau$$
Let's write $\omega_\tau(x) = |x|^{i\tau}$. Then
$$F(x) = -{1 \over {2\pi}}\int_{-\infty}^{+\infty} \widehat{g}({1\over2} + i\tau)\Lambda(\chi\omega_\tau, {1\over2})\,\overline{\chi\omega_\tau(x)}\, d\tau$$
We now split the discussion according to whether $\nu$ is finite or archimedean.

If $\nu$ is finite, we see that the integrand has period $2\pi\over\log(q)$ except for $\widehat{g}$. The Poisson summation formula proven before implies
$$\eqalign{
\sum_{j\in\ZZ}\widehat{g}(j\,{2\pi i\over\log(q)} + {1\over2} + i\tau)
&= \log(q)\,\sum_{k\in\ZZ}g(q^k)\,q^{k({1\over2} + i\tau)} \cr
&= \log(q)\, \int_{K_\nu^\times}g(|x|) |x|^{{1\over2} + i\tau}\,d^*x \cr
&= \log(q)\ f(\chi\omega_\tau) }$$
so that we obtain in the end
$$F(x) = -{\log(q) \over {2\pi}}\int_{-{\pi\over\log(q)}}^{+{\pi\over\log(q)}} f(\chi\omega_\tau)\Lambda(\chi\omega_\tau, {1\over2})\,\overline{\chi\omega_\tau(x)}\, d\tau$$
Starting from
$$g(|x|) = {1 \over {2\pi i}}\int_{Re(s)=1/2} \widehat{g}(s)|x|^{-s}\, ds$$
leads in the same manner to
$$f(x) = {\log(q) \over {2\pi}}\int_{-{\pi\over\log(q)}}^{+{\pi\over\log(q)}} f(\chi\omega_\tau)\,\overline{\chi\omega_\tau(x)}\, d\tau$$
This identifies in a concrete way that component of $S_\nu$ containing $\chi(1/x)$ and proves the Theorem in the finite case.

In the archimedean case we can just identify immediately $\widehat{g}({1\over2} + i\tau)$ with $f(\chi\omega_\tau)$ as the multiplicative Haar measure $d^*x$ has been chosen to be pushed forward to $du/u$ under the norm map $x\mapsto|x|$. So the formulae in that case are
$$F(x) = -{1 \over {2\pi}}\int_{-\infty}^{+\infty} f(\chi\omega_\tau)\Lambda(\chi\omega_\tau, {1\over2})\,\overline{\chi\omega_\tau(x)}\, d\tau$$
$$f(x) = {1\over {2\pi}}\int_{-\infty}^{+\infty} f(\chi\omega_\tau)\,\overline{\chi\omega_\tau(x)}\, d\tau$$
This completes the proof of the Theorem.

\par\vfill\eject
{
{\bf REFERENCES}\par
\baselineskip = 12 pt
\parskip = 4 pt
\font\smallRoman = cmr8
\smallRoman
\font\smallBold = cmbx8
\font\smallSlanted = cmsl8
{\smallBold [1] J. Tate}, Thesis, Princeton 1950, reprinted in Algebraic Number Theory, ed. J.W.S. Cassels and A. Fr\"ohlich, Academic Press, (1967).\par
{\smallBold [2] A. Weil},{\smallSlanted ``Sur les ``formules explicites'' de la th\'eorie des nombres premiers''}, Comm. Lund (vol d\'edi\'e \`a Marcel Riesz), (1952).\par
{\smallBold [3] S. Haran},{\smallSlanted ``Riesz potentials and explicit sums in arithmetic''}, Invent. Math.  101, 697-703 (1990).\par
{\smallBold [4] J.F. Burnol}, {\smallSlanted ``The Explicit Formula and a Propagator''}, electronic manuscript, available at the http://xxx.lanl.gov server, math/9809119 (September 1998, revised November 1998)\par
{\smallBold [5] I. M. Gel'fand, M. I. Graev, I. I. Piateskii-Shapiro},{\smallSlanted ``Representation Theory and automorphic functions''}, Philadelphia, Saunders (1969).\par
{\smallBold [6] A. Weil},{\smallSlanted ``Fonctions z\^etas et distributions''}, S\'eminaire Bourbaki n${}^{\smallRoman o}$ 312, (1966).\par
{\smallBold [7] V. S. Vladimirov},{\smallSlanted ``Generalized functions over the fields of p-adic numbers''}, Russian Math. Surveys 43:5, 19-64 (1988).\par
\vfill
\centerline{Jean-Fran\c{c}ois Burnol}
\centerline{62 rue Albert Joly}
\centerline{F-78000 Versailles}
\centerline{France}
\centerline{jf.burnol@dial.oleane.com}
\centerline{November 1998 (revised)}
}
\eject
\bye